\documentclass{amsart}

\usepackage{amssymb}
\usepackage{amsmath}
\usepackage{amsthm}
\usepackage{graphics}

\allowdisplaybreaks

\numberwithin{equation}{section}

\theoremstyle{plain}
\newtheorem{thm}{Theorem}[section]
\newtheorem{prop}[thm]{Proposition}
\newtheorem{lem}[thm]{Lemma}

\theoremstyle{definition}

\newcommand{\R}{\mathbb{R}}

\newcommand{\Z}{\mathbb{Z}}

\newcommand{\calB}{\mathcal{B}}
\newcommand{\calD}{\mathcal{D}}
\newcommand{\calF}{\mathcal{F}}
\newcommand{\calL}{\mathcal{L}}
\newcommand{\calM}{\mathcal{M}}
\newcommand{\calS}{\mathcal{S}}

\begin{document}

\title[Pseudo-differential operators on modulation spaces]
{Boundedness properties of pseudo-differential operators and
Calder\'on-Zygmund operators \\ on modulation spaces}
\author{Mitsuru Sugimoto \and Naohito Tomita}
\date{}

\address{Mitsuru Sugimoto \\
Department of Mathematics \\
Graduate School of Science \\
Osaka University \\
Toyonaka, Osaka 560-0043, Japan}
\email{sugimoto@math.wani.osaka-u.ac.jp}

\address{Naohito Tomita \\
Department of Mathematics \\
Graduate School of Science \\
Osaka University \\
Toyonaka, Osaka 560-0043, Japan}
\email{tomita@gaia.math.wani.osaka-u.ac.jp}

\keywords{Calder\'on-Zygmund operators,
modulation spaces, pseudo-differential operators}

\subjclass[2000]{42B20, 42B35, 47G30}

\begin{abstract}
In this paper,
we study the boundedness of pseudo-differential operators
with symbols in $S_{\rho,\delta}^m$ on the modulation spaces $M^{p,q}$.
We discuss the order $m$ for the boundedness
$\mathrm{Op}(S_{\rho,\delta}^m) \subset \calL(M^{p,q}(\R^n))$
to be true.
We also prove the existence of a Calder\'on-Zygmund operator
which is not bounded
on the modulation space $M^{p,q}$ with $q \neq 2$.
This unboundedness is still true even if we assume
a generalized $T(1)$ condition.
These results are induced by the unboundedness of
pseudo-differential operators on $M^{p,q}$
whose symbols are of the class $S_{1,\delta}^0$ with $0<\delta<1$.
\end{abstract}

\maketitle

\section{Introduction}\label{section1}
The modulation spaces $M^{p,q}$ introduced by Feichtinger
\cite{Feichtinger1,Feichtinger2} are
fundamental function spaces of time-frequency
analysis which is originated in signal analysis or quantum mechanics.
See also \cite{Feichtinger3} or Triebel \cite{Triebel}.
Recently these spaces have been also recognized as a
useful tool for the theory of pseudo-differential operators
(see Gr\"ochenig \cite{Grochenig2}).
The objective of this paper is to discuss the boundedness of
pseudo-differential operators and also
the unboundedness of Calder\'on-Zygmund operators on modulation spaces.
\par
The boundedness of pseudo-differential operators on the modulation spaces
was studied by many authors, for example,
Gr\"ochenig and Heil \cite{Grochenig-Heil},
Tachizawa \cite{Tachizawa}, and Toft \cite{Toft}.
They proved that pseudo-differential operators
with symbols in some modulation space are $M^{p,q}$-bounded,
and as a corollary we have the $M^{p,q}$-boundedness of
pseudo-differential operators with symbols in H\"ormander's class $S_{0,0}^0$.
A pioneering work of Sj\"ostrand \cite{Sjostrand} should be also
mentioned here which proved the $L^2$-boundedness by introducing
a symbol class based on the spirit of time-frequency analysis.
We note here that $L^2=M^{2,2}$.
\par
In this paper, we consider the case of general symbol classes
$S_{\rho,\delta}^m$.
First we recall the result of Calder\'on and Vaillancourt \cite{C-V}, which
shows that pseudo-differential operators
with symbols in $S_{\delta,\delta}^0$ with $0<\delta<1$
(hence in $S_{\rho,\delta}^0$ with $0\leq\delta\leq\rho\leq1$, $\delta<1$)
are $L^2$-bounded
by reducing them to the case of $S_{0,0}^0$.
The proof was carried out by \lq\lq dilation argument'' based on the fact
\[
\|f(\lambda x)\|_{L^2(\R)}=\lambda^{-n/2}\|f(x)\|_{L^2(\R)},\qquad \lambda>0.
\]
We can conclude that the operator norm of $\sigma(X,D)$ is equal to that of
$\sigma(\lambda X,\lambda^{-1}D)$ from this equality.
By following the same argument, we can also expect that pseudo-differential
operators with symbols in $S_{\delta,\delta}^0$ with $0<\delta<1$
are $M^{p,q}$-bounded.
However, the dilation property of the modulation spaces $M^{p,q}$
is rather different as was investigated in author's previous paper
\cite{Sugimoto-Tomita-1}, and due to the property, we can give the
following negative answer to this expectation:
\begin{thm}\label{1.1}
Let $1<q<\infty$,
$m \in \R$,
$0 \le \delta \le \rho \le 1$ and $\delta<1$.
Then we have the boundedness
$\mathrm{Op}(S_{\rho,\delta}^m) \subset \calL(M^{2,q}(\R^n))$
if and only if $m \le -|1/q-1/2|\delta n$.
\end{thm}
The $\lq\lq$only if'' part of Theorem \ref{1.1}
is a restricted version of Theorem \ref{3.6} 
which says that
$\mathrm{Op}(S_{\rho,\delta}^m) \subset \calL(M^{p,q}(\R^n))$
only if $m \le -|1/q-1/2|\delta n$.
This result was first proved in author's paper \cite{Sugimoto-Tomita-2}.
The $\lq\lq$if" part is a restricted version
of Theorem \ref{3.1} in Section \ref{section3}
which treats the general $M^{p,q}$-boundedness
with $m\leq m(p,q)$, where $m(p,q)$ is an index such that
$m(2,q)=-|1/q-1/2|\delta n$.
We remark that the boundedness
$\mathrm{Op}(S_{\rho,\delta}^m) \subset \calL(M^{p,q}(\R^n))$
with $m < -|1/q-1/2|\delta n$ (Proposition \ref{3.2})
is a straightforward consequence of the dilation argument stated in the above
if we directly use the dilation property
(Proposition \ref{2.1}) of the modulation spaces.
The main contribution of Theorem \ref{1.1}
is the boundedness result with the 
critical order for $m$.
\par
Theorem \ref{1.1} also indicates
a difference between $L^p$ spaces and the modulation spaces $M^{p,q}$.
Fefferman \cite{Fefferman} proved that
$\mathrm{Op}(S_{\rho,\delta}^m) \subset \calL(L^p(\R^n))$
if $m \le -|1/p-1/2|(1-\rho)n$, $0 \le \delta \le \rho \le 1$
and $\delta<1$.
Moreover, it is known that this order of $m$ is critical
(\cite[Chapter 7, 5.12]{Stein}).
We remark that,
for the $L^p$-boundedness
of pseudo-differential operators with symbols in $S_{\rho,\delta}^m$,
the critical order of $m$ is determined by $\rho$.
On the other hand, Theorem \ref{1.1}
says that,
for the $M^{p,q}$-boundedness,
the critical order of $m$ is determined by $\delta$
(at least, in the case $p=2$).
\medskip
\par
By developing the investigation of pseudo-differential operators,
we can also discuss the boundedness property of
Calder\'on-Zygmund operators on the modulation spaces. 
A Calder\'on-Zygmund operator is an $L^2$-bounded linear mapping
whose distributional kernel is
a function $K(x,y)$ outside the diagonal $\{x=y\}$ satisfying
\begin{align*}
&|K(x,y)| \le
C|x-y|^{-n},
\\
&|K(x,y)-K(x',y)|
\le C\frac{|x-x'|^{\epsilon}}{|x-y|^{n+\epsilon}}
\qquad\text{for $|x-y| >2|x-x'|$}, 
\\
&|K(x,y)-K(x,y')|
\le C\frac{|y-y'|^{\epsilon}}{|x-y|^{n+\epsilon}}
\qquad\text{for $|x-y| >2|y-y'|$},
\end{align*}
where $0<\epsilon\le 1$.
It is well know that Calder\'on-Zygmund operators
are $L^p$-bounded for $1<p<\infty$
(see, for example, \cite[Theorem 5.10]{Duo}).
In this paper, we consider the problem
\lq\lq Are Calder\'on-Zygmund operators $M^{p,q}$-bounded?".
It is easy to show that the $L^p$-boundedness of
the operators of convolution type induces the $M^{p,q}$-boundedness.
Hence, Calder\'on-Zygmund operators of convolution type
(Riesz transforms, for example) are always $M^{p,q}$-bounded.
We will prove that it is not generally true
if remove the convolution type assumption.
The following theorem is a simplified version
of Theorem \ref{4.1} in Section \ref{section4}.
\begin{thm}\label{1.2}
Let $1<p,q<\infty$.
If $q \neq 2$, then there exists a Calder\'on-Zygmund operator
which is not bounded on $M^{p,q}(\R^n)$.
\end{thm}
Theorem \ref{1.2} is a direct consequence of the unboundedness of
pseudo-differential operators on the modulation space $M^{p,q}$
whose symbols are of the class $S_{1,\delta}^0$ with $0<\delta<1$
(see Theorem \ref{3.6}).
In fact, it is know that such pseudo-differential
operators are Calder\'on-Zygmund operators (see
\cite[Chapter 7, Sections 1 and 2]{Stein}).

On the other hand, we know some boundedness results on the Besov
spaces $\Dot{B}^{s,q}_p$ and Triebel-Lizorkin spaces $\Dot{F}^{s,q}_p$.
By the same reason described above,
Calder\'on-Zygmund operators of convolution type
are always $\Dot{B}^{s,q}_p$-bounded.
Furthermore, in the case of non-convolution type,
the following generalized $T(1)$ condition
is useful to discuss the $\Dot{B}^{s,q}_p$-boundedness of $T$:
\begin{equation}\label{(1.1)}
T(P)=0\qquad\text{for all polynomials $P$ such that
$\operatorname{deg}P\leq \ell$}.
\end{equation}
This condition, together with an extra condition on the smoothness
of the kernel of $T$
and the {\it weak boundedness property}
which will be defined in Section \ref{section4},
induce the $\Dot{B}^{s,q}_p$-boundedness of $T$,
where $s$ is determined by the order of the polynomials $\ell$
and the smoothness order of the kernel
(see Lemari\'e \cite{Lemarie},
Meyer and Coifman \cite[p.114]{M-C}).
For the $\Dot{F}^{s,q}_p$-boundedness, we need more conditions
on the transpose $T^*$ of $T$, that is, a smoothness condition
of the kernel of $T^*$ and the condition
\begin{equation}\label{(1.2)}
T^*(P)=0\qquad\text{for all polynomials $P$ such that
$\operatorname{deg}P\leq \ell^*$}
\end{equation}
(see Frazier, Torres and Weiss \cite{F-T-W}).
We remark that the boundedness on the homogeneous spaces
($\Dot{B}^{s,q}_p, \Dot{F}^{s,q}_p$)
induces that on the inhomogeneous spaces
($B^{s,q}_p, F^{s,q}_p$)
under suitable conditions.

On account of these results, we can expect the $M^{p,q}$-boundedness
of Calder\'on-Zygmund operators $T$ which satisfy
the smoothness conditions on the kernels,
conditions \eqref{(1.1)}-\eqref{(1.2)},
and the weak boundedness property.
But even if we assume these reasonable conditions,
we can never prove the $M^{p,q}$-boundedness of
Calder\'on-Zygmund operators.
See Theorem \ref{4.1} for the detailed statement.

\section{Preliminaries}\label{section2}
Let $\calS(\R^n)$ and $\calS'(\R^n)$ be the Schwartz spaces of
all rapidly decreasing smooth functions
and tempered distributions,
respectively.
We define the Fourier transform $\calF f$
and the inverse Fourier transform $\calF^{-1}f$
of $f \in \calS(\R^n)$ by
\[
\calF f(\xi)
=\widehat{f}(\xi)
=\int_{\R^n}e^{-i\xi \cdot x}\, f(x)\, dx
\quad \text{and} \quad
\calF^{-1}f(x)
=\frac{1}{(2\pi)^n}
\int_{\R^n}e^{ix\cdot \xi}\, f(\xi)\, d\xi.
\]

We introduce the modulation spaces
based on Gr\"ochenig \cite{Grochenig}.
Fix a function $\gamma \in \calS(\R^n)\setminus \{ 0 \}$
(called the {\it window function}).
Then the short-time Fourier transform $V_{\gamma}f$ of
$f \in \calS'(\R^n)$ with respect to $\gamma$
is defined by
\[
V_{\gamma}f(x,\xi)
=(f, M_{\xi}T_x \gamma)
\qquad \text{for} \ x, \xi \in \R^n,
\]
where $M_{\xi}\gamma(t)=e^{i\xi \cdot t}\gamma(t)$,
$T_x \gamma(t)=\gamma(t-x)$
and $(\cdot,\cdot)$ denotes the inner product on $L^2(\R^n)$.
We note that,
for $f \in \calS'(\R^n)$,
$V_{\gamma}f$ is continuous on $\R^{2n}$
and $|V_{\gamma}f(x,\xi)|\le C(1+|x|+|\xi|)^N$
for some constants $C,N \ge 0$
(\cite[Theorem 11.2.3]{Grochenig}).
Let $1\le p,q \le \infty$.
Then the modulation space $M^{p,q}(\R^n)$
consists of all $f \in \calS'(\R^n)$
such that
\[
\|f\|_{M^{p,q}}
=\|V_{\gamma}f\|_{L^{p,q}}
=\left\{ \int_{\R^n} \left(
\int_{\R^n} |V_{\gamma}f(x,\xi)|^{p}\, dx
\right)^{q/p} d\xi \right\}^{1/q}
< \infty.
\]
We note that
$M^{2,2}(\R^n)=L^2(\R^n)$
(\cite[Proposition 11.3.1]{Grochenig})
and $M^{p,q}(\R^n)$ is a Banach space
(\cite[Proposition 11.3.5]{Grochenig}).
The definition of $M^{p,q}(\R^n)$ is independent
of the choice of the window function
$\gamma \in \calS(\R^n)\setminus \{ 0 \}$,
that is,
different window functions
yield equivalent norms
(\cite[Proposition 11.3.2]{Grochenig}).
We denote by $\calL(M^{p,q}(\R^n))$ the space of
all bounded linear operators on $M^{p,q}(\R^n)$.

In order to state the dilation property of the modulation spaces,
we introduce the indices.
For $1\le p \le \infty$,
$p'$ is the conjugate exponent of $p$
(that is, $1/p+1/p'=1$).
We define subsets of $(1/p,1/q)\in[0,1]\times[0,1]$ in the following way:
\begin{alignat*}{2}
&I_1\,:\,\min(1/q,1/2)\ge 1/p, &\qquad
&I_1^*\,:\,\max(1/q,1/2)\le 1/p,
\\
&I_2\,:\,\min(1/p,1/p')\ge 1/q, &\qquad
&I_2^*\,:\,\max(1/p,1/p')\le 1/q, 
\\
&I_3\,:\,\min(1/q,1/2)\ge 1/p', &\qquad
&I_3^*\,:\,\max(1/q,1/2)\le 1/p'.
\end{alignat*}
See the following figure:
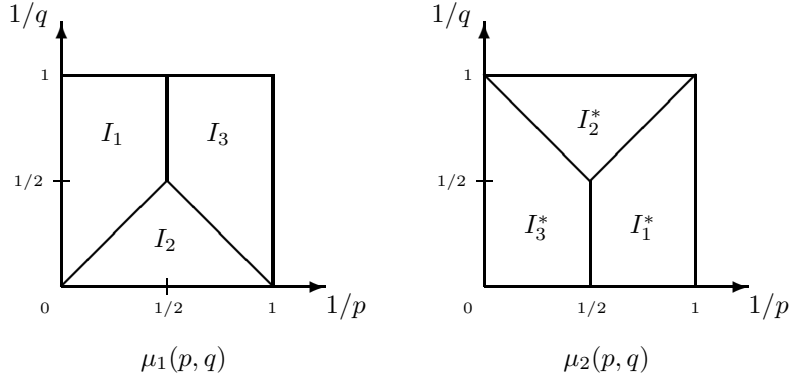
\begin{center}
\begin{picture}(300,160)
\thicklines
\put(20,40){\vector(1,0){100}}
\put(20,120){\line(1,0){80}}
\put(20,40){\vector(0,1){100}}
\put(100,40){\line(0,1){80}}
\put(60,80){\line(0,1){40}}
\put(20,40){\line(1,1){40}}
\put(60,80){\line(1,-1){40}}
\put(50,10){$\mu_1(p,q)$}
\put(120,30){$1/p$}
\put(0,140){$1/q$}
\put(35,95){$I_1$}
\put(55,55){$I_2$}
\put(75,95){$I_3$}
\thinlines
\put(17,80){\line(1,0){6}}
\put(60,37){\line(0,1){6}}
\put(12,30){{\tiny $0$}}
\put(55,30){{\tiny $1/2$}}
\put(98,30){{\tiny $1$}}
\put(2,78){{\tiny $1/2$}}
\put(12,118){{\tiny $1$}}
\thicklines
\put(180,40){\vector(1,0){100}}
\put(180,120){\line(1,0){80}}
\put(180,40){\vector(0,1){100}}
\put(260,40){\line(0,1){80}}
\put(220,40){\line(0,1){40}}
\put(180,120){\line(1,-1){40}}
\put(220,80){\line(1,1){40}}
\put(210,10){$\mu_2(p,q)$}
\put(280,30){$1/p$}
\put(160,140){$1/q$}
\put(235,60){$I_1^*$}
\put(215,100){$I_2^*$}
\put(195,60){$I_3^*$}
\thinlines
\put(177,80){\line(1,0){6}}
\put(172,30){{\tiny $0$}}
\put(215,30){{\tiny $1/2$}}
\put(258,30){{\tiny $1$}}
\put(162,78){{\tiny $1/2$}}
\put(172,118){{\tiny $1$}}
\end{picture}
\end{center}
We introduce the indices
\begin{align*}
&\mu_1(p,q)
=\max\{0,1/q-\min(1/p,1/p')\}-1/p,
\\
&\mu_2(p,q)
=\min\{0,1/q-\max(1/p,1/p')\}-1/p.
\end{align*}
Then
\begin{align*}
&\mu_1(p,q)=
\begin{cases}
-2/p+1/q &\text{if} \quad(1/p,1/q) \in I_1, \\
-1/p &\text{if} \quad(1/p,1/q) \in I_2, \\
1/q-1 &\text{if} \quad(1/p,1/q) \in I_3,
\end{cases}
\\
&\mu_2(p,q)=
\begin{cases}
-2/p+1/q &\text{if} \quad(1/p,1/q) \in I_1^*, \\
-1/p &\text{if} \quad(1/p,1/q) \in I_2^*, \\
1/q-1 &\text{if} \quad(1/p,1/q) \in I_3^*.
\end{cases}
\end{align*}
We define the dilation operator $\Lambda_a$ by $\Lambda_a f(x)=f(a x)$,
where $a>0$.
The following proposition plays an important role
in the proofs of Theorem \ref{3.1} and Proposition \ref{3.2}.
\begin{prop}[{\cite[Theorem 1.1]{Sugimoto-Tomita-1}}]\label{2.1}
Let $1\le p,q \le \infty$.
Then the following are true:
\begin{enumerate}
\item
There exists a constant $C>0$ such that
\[
\|\Lambda_a f\|_{M^{p,q}}
\le Ca^{n\mu_1(p,q)}\|f\|_{M^{p,q}}
\quad \text{for all} \
f \in M^{p,q}(\R^n)
\ \text{and} \
a \ge 1.
\]
\item
There exists a constant $C>0$ such that
\[
\|\Lambda_a f\|_{M^{p,q}}
\le Ca^{n\mu_2(p,q)}\|f\|_{M^{p,q}}
\quad \text{for all} \
f \in M^{p,q}(\R^n)
\ \text{and} \
0<a \le 1.
\]
\end{enumerate}
\end{prop}

The optimality of the power of $a$ in Proposition \ref{2.1}
is also discussed in \cite{Sugimoto-Tomita-1}.

\section{The boundedness of pseudo-differential operators}\label{section3}
Let $m \in \R$ and $0 \le \delta \le \rho \le 1$.
The symbol class $S_{\rho,\delta}^m$ consists of
all $\sigma \in C^{\infty}(\R^n \times \R^n)$ such that
\[
|\partial_{\xi}^{\alpha}\partial_x^{\beta}\sigma(x,\xi)|
\le C_{\alpha,\beta}(1+|\xi|)^{m-\rho |\alpha|+\delta |\beta|}
\]
for all $\alpha,\beta \in \Z_+^n=\{0,1,\dots\}^n$.
We denote by
$|\cdot|_{S_{\rho,\delta}^m,N}, \ N=0,1,\dots$,
the semi-norms on $S_{\rho,\delta}^m$,
that is,
\[
|\sigma|_{S_{\rho,\delta}^m,N}
=\max_{|\alpha+\beta| \le N}
\sup_{x,\xi \in \R^n}(1+|\xi|)^{-(m-\rho|\alpha|+\delta|\beta|)}
|\partial_{\xi}^{\alpha}\partial_x^{\beta}\sigma(x,\xi)|.
\]
For $\sigma \in S_{\rho,\delta}^m$,
the pseudo-differential operator $\sigma(X,D)$ is defined by
\[
\sigma(X,D)f(x)
=\frac{1}{(2\pi)^n}
\int_{\R^n}e^{ix\cdot \xi}\,
\sigma(x,\xi)\, \widehat{f}(\xi)\, d\xi
\]
for $f \in \calS(\R^n)$.
We denote by $\mathrm{Op}(S_{\rho,\delta}^m)$
the class of all pseudo-differential operators
with symbols in $S_{\rho,\delta}^m$.
Given a symbol $\sigma \in S_{\rho,\delta}^m$
with $\delta<1$,
the symbol $\sigma^*$ defined by
\begin{align}\label{(3.1)}
\sigma^*(x,\xi)
&=\textrm{Os-}
\frac{1}{(2\pi)^n}\int_{\R^n}\int_{\R^n}
e^{-iy\cdot \zeta}\,
\overline{\sigma(x+y,\xi+\zeta)}\,
dy \, d\zeta \\
&=\lim_{\epsilon \to 0}
\frac{1}{(2\pi)^n}
\int_{\R^n}\int_{\R^n}
e^{-iy\cdot \zeta}\,
\chi(\epsilon y, \epsilon \zeta)\,
\overline{\sigma(x+y,\xi+\zeta)}\,
dy \, d\zeta \nonumber
\end{align}
satisfies $\sigma^* \in S_{\rho,\delta}^m$
and
\begin{equation}\label{(3.2)}
(\sigma(X,D)f,g)=(f,\sigma^*(X,D)g)
\qquad \text{for all} \ f,g \in \calS(\R^n),
\end{equation}
where $\chi \in \calS(\R^{2n})$
satisfies $\chi(0,0)=1$
(\cite[Chapter 2, Theorem 2.6]{Kumano-go}).
Note that oscillatory integrals
are independent of the choice of $\chi \in \calS(\R^{2n})$
satisfying $\chi(0,0)=1$
(\cite[Chapter 1, Theorem 6.4]{Kumano-go}),
and the derivatives of $\sigma^*(x,\xi)$ can be written as
\begin{equation}\label{(3.3)}
\partial_{\xi}^{\alpha}\partial_x^{\beta}\sigma^*(x,\xi)
=\textrm{Os-}
\frac{1}{(2\pi)^n}\int_{\R^n}\int_{\R^n}
e^{-iy\cdot \zeta}\,
\overline{(\partial_{\xi}^{\alpha}\partial_x^{\beta}\sigma)(x+y,\xi+\zeta)}\,
dy \, d\zeta
\end{equation}
in virtue of
\cite[Chapter 1, Theorem 6.6]{Kumano-go}
(see also \cite[p.70, (2.23)]{Kumano-go}).
\par
Our main result on the boundedness of pseudo-differential operators
is the following:
\begin{thm}\label{3.1}
Let $1<p,q<\infty$,
$m \in \R$,
$0 \le \delta \le \rho \le 1$ and $\delta<1$.
If $m \le -(\mu_1(p,q)-\mu_2(p,q))\delta n$,
then $\mathrm{Op}(S_{\rho,\delta}^m) \subset \calL(M^{p,q}(\R^n))$.
\end{thm}
In order to clarify $\mu_1(p,q)-\mu_2(p,q)$,
we divide $I_1,\dots,I_3^*$ in the following way:
\[
J_1\,:\, (I_1 \cap I_2^*) \cup (I_2 \cap I_1^*), \quad
J_2\,:\, (I_1 \cap I_3^*) \cup (I_3 \cap I_1^*), \quad
J_3\,:\, (I_2 \cap I_3^*) \cup (I_3 \cap I_2^*).
\]
See the following figure:
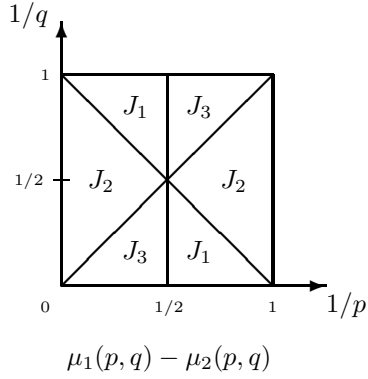
\begin{center}
\begin{picture}(150,160)
\thicklines
\put(20,40){\vector(1,0){100}}
\put(20,120){\line(1,0){80}}
\put(20,40){\vector(0,1){100}}
\put(100,40){\line(0,1){80}}
\put(60,40){\line(0,1){80}}
\put(20,40){\line(1,1){80}}
\put(20,120){\line(1,-1){80}}
\put(22,10){$\mu_1(p,q)-\mu_2(p,q)$}
\put(120,30){$1/p$}
\put(0,140){$1/q$}
\put(43,105){$J_1$}
\put(30,77){$J_2$}
\put(43,50){$J_3$}
\put(67,50){$J_1$}
\put(80,77){$J_2$}
\put(67,105){$J_3$}
\thinlines
\put(17,80){\line(1,0){6}}
\put(12,30){{\tiny $0$}}
\put(55,30){{\tiny $1/2$}}
\put(98,30){{\tiny $1$}}
\put(2,78){{\tiny $1/2$}}
\put(12,118){{\tiny $1$}}
\end{picture}
\end{center}
Then we have
\begin{align*}
&\mu_1(p,q)-\mu_2(p,q)=
\begin{cases}
|1/p-1/q| &\text{if} \quad(1/p,1/q) \in J_1, \\
|2/p-1| &\text{if} \quad(1/p,1/q) \in J_2, \\
|1/p+1/q-1| &\text{if} \quad(1/p,1/q) \in J_3.
\end{cases}
\end{align*}
Let $\psi_0 ,\psi \in \calS(\R^n)$ be such that
$\mathrm{supp}\, \psi_0 \subset \{\xi: |\xi|\le 2\}$,
$\mathrm{supp}\, \psi \subset \{\xi: 1/2 \le |\xi| \le 2\}$
and $\psi_0(\xi)+\sum_{j=1}^{\infty}\psi(2^{-j}\xi)=1$
for all $\xi \in \R^n$.
Set $\psi_j=\psi(2^{-j}\cdot)$ if $j \ge 1$.
Then
\begin{equation}\label{(3.4)}
\|\partial^{\alpha}\psi_j\|_{L^{\infty}}
\le C_{\alpha}2^{-j|\alpha|},
\quad \mathrm{supp}\, \psi_j \subset \{2^{j-1}\le |\xi|\le 2^{j+1}\},
\quad \sum_{j=0}^{\infty}\psi_j \equiv 1.
\end{equation}
The following weak form of Theorem \ref{3.1} (when $\delta<1$)
is a straightforward consequence of the dilation property
of the modulation spaces (Proposition \ref{2.1}).
\begin{prop}\label{3.2}
Let $1<p,q<\infty$,
$m \in \R$ and $0 \le \delta \le \rho \le 1$.
If $m < -(\mu_1(p,q)-\mu_2(p,q))\delta n$,
then $\mathrm{Op}(S_{\rho,\delta}^m) \subset \calL(M^{p,q}(\R^n))$.
\end{prop}
\begin{proof}
Our proof is based on that of \cite[Chapter 7, Theorem 2]{Stein}.
Let $\sigma \in S_{\rho,\delta}^m$.
By the decomposition \eqref{(3.4)},
we have
\[
\sigma(x,\xi)=\sum_{j=0}^{\infty}\sigma(x,\xi)\, \psi_j(\xi)
=\sum_{j=0}^{\infty}\sigma_j(x,\xi).
\]
Since
$\sigma_j(X,D)
=\Lambda_{2^{j\delta}}\, \sigma_j(2^{-j\delta}X, 2^{j\delta}D)\,
\Lambda_{2^{-j\delta}}$,
by Proposition \ref{2.1},
we see that
\begin{align*}
\|\sigma(X,D)f\|_{M^{p,q}}
&\le \sum_{j=0}^{\infty}\|\sigma_j(X,D)f\|_{M^{p,q}} \\
&=\sum_{j=0}^{\infty}
\|\Lambda_{2^{j\delta}}\, \sigma_j(2^{-j\delta}X, 2^{j\delta}D)\,
\Lambda_{2^{-j\delta}}f\|_{M^{p,q}} \\
&\le \sum_{j=0}^{\infty}2^{j\delta n(\mu_1(p,q)-\mu_2(p,q))}
\|\sigma_j(2^{-j\delta}X, 2^{j\delta}D)\|_{\calL(M^{p,q})}
\|f\|_{M^{p,q}}
\end{align*}
for all $f \in \calS(\R^n)$.
Since $S_{\rho,\delta}^m \subset S_{\delta,\delta}^m$
and $1+2^{j\delta}|\xi| \sim 2^j$
on $\mathrm{supp}\, \psi_j(2^{j\delta}\cdot)$,
if we set
$\tau_j(x,\xi)=2^{-jm}\sigma_j(2^{-j\delta}x, 2^{j\delta}\xi)$,
then for any $\alpha,\beta \in \Z_+^n$
\[
|\partial_{\xi}^{\alpha}\partial_x^{\beta}\tau_j(x,\xi)|
\le C_{\alpha,\beta}|\sigma|_{S_{\rho,\delta}^m, |\alpha+\beta|}
\]
for all $j \in \Z_+$.
Hence,
by $\mathrm{Op}(S_{0,0}^0) \subset \calL(M^{p,q}(\R^n))$
and $m+(\mu_1(p,q)-\mu_2(p,q))\delta n<0$,
we have
\begin{align*}
&\sum_{j=0}^{\infty}2^{j\delta n(\mu_1(p,q)-\mu_2(p,q))}
\|\sigma_j(2^{-j\delta}X, 2^{j\delta}D)\|_{\calL(M^{p,q})}
\|f\|_{M^{p,q}} \\
&=\sum_{j=0}^{\infty}2^{j(m+(\mu_1(p,q)-\mu_2(p,q))\delta n)}
\|\tau_j(X,D)\|_{\calL(M^{p,q})}
\|f\|_{M^{p,q}} \\
&\le C
\left( \sum_{j=0}^{\infty}2^{j(m+(\mu_1(p,q)-\mu_2(p,q))\delta n)}\right)
\|f\|_{M^{p,q}}
=C\|f\|_{M^{p,q}}
\end{align*}
for all $f \in \calS(\R^n)$.
The proof is complete.
\end{proof}
In view of Proposition \ref{3.2},
the non-trivial part of Theorem \ref{3.1}
is the boundedness with the critical order
$m=-(\mu_1(p,q)-\mu_2(p,q))\delta n$.
The rest of this section is devoted to
the proof of it.
Let $\varphi \in \calS(\R^n)$ be such that
\begin{equation}\label{(3.5)}
\mathrm{supp}\, \varphi \subset [-1,1]^n
\quad \text{and} \quad
\sum_{k \in \Z^n}\varphi(\xi-k) \equiv 1.
\end{equation}
By the decompositions \eqref{(3.4)} and \eqref{(3.5)},
we have
\[
1 \equiv \sum_{j=0}^{\infty}\psi_j(\xi)
=\sum_{j=0}^{\infty}\psi_j(\xi)
\left(\sum_{k \in \Z^n}\varphi(2^{-j\delta}y-k)\right)
\left(\sum_{\ell \in \Z^n}\varphi(2^{-j\delta}\xi-\ell)\right).
\]
Hence,
\begin{equation}\label{(3.6)}
\sum_{k \in \Z^n}\sum_{\ell \in \Z^n}\sum_{j=0}^{\infty}
\varphi(2^{-j\delta}y-k)\,
\varphi(2^{-j\delta}\xi-\ell)\, \psi_j(\xi)=1
\end{equation}
for all $(y,\xi) \in \R^n \times \R^n$.
For $\sigma \in S_{\rho,\delta}^m$,
we set
\begin{equation}\label{(3.7)}
\sigma_j^{k,\ell}(x,\xi)
=\varphi(2^{-j\delta}D_x-k)\, \sigma(x,\xi)\,
\varphi(2^{-j\delta}\xi-\ell)\, \psi_j(\xi),
\end{equation}
where
\begin{align*}
\varphi(2^{-j\delta}D_x-k)\, \sigma(x,\xi)
&=\calF_1^{-1}[\varphi(2^{-j\delta}\cdot-k)\,
\calF_1\sigma(\cdot,\xi)](x) \\
&=\frac{1}{(2\pi)^n}
\int_{\R^n}e^{ix\cdot y}\,
\varphi(2^{-j\delta}y-k)\,
\calF_1\sigma(y,\xi)\, dy,
\end{align*}
$\calF_1$ and $\calF_1^{-1}$ are the partial Fourier transform 
and inverse Fourier transform in the first variable,
respectively.
\begin{lem}\label{3.3}
Let $1 \le p \le \infty$,
$m \in \R$,
$0 \le \delta \le \rho \le 1$
and $\sigma \in S_{\rho,\delta}^m$.
Then there exists a constant $C>0$ such that
\[
\|\sigma_j^{k,\ell}(X,D)f\|_{L^p}
\le C2^{jm}(1+|k|)^{-n-1}\|f\|_{L^p}
\]
for all $j \in \Z_+$, $k, \ell \in \Z^n$ and
$f \in \calS(\R^n)$,
where $\sigma_j^{k,\ell}$ is defined by \eqref{(3.7)}.
\end{lem}
\begin{proof}
Using
$\sigma_j^{k,\ell}(X,D)
=\Lambda_{2^{j\delta}}\,
\sigma_j^{k,\ell}(2^{-j\delta}X,2^{j\delta}D)\,
\Lambda_{2^{-j\delta}}$,
we have
\[
\|\sigma_j^{k,\ell}(X,D)f\|_{L^p}
=2^{-j\delta n/p}
\|\sigma_j^{k,\ell}(2^{-j\delta}X,2^{j\delta}D)
\Lambda_{2^{-j\delta}}f\|_{L^p}
\]
for all $f \in \calS(\R^n)$.
Set
\[
K_j^{k,\ell}(x,y)
=\frac{1}{(2\pi)^n}\int_{\R^n}e^{i(x-y)\cdot\xi}\,
\sigma_j^{k,\ell}(2^{-j\delta}x,2^{j\delta}\xi)\, d\xi.
\]
Since
\[
\sigma_j^{k,\ell}(2^{-j\delta}X,2^{j\delta}D)\Lambda_{2^{-j\delta}}f(x)
=\int_{\R^n}K_j^{k,\ell}(x,y)\, (\Lambda_{2^{-j\delta}}f)(y)\, dy,
\]
if
\begin{equation}\label{(3.8)}
|K_j^{k,\ell}(x,y)|
\le C2^{jm}(1+|k|)^{-n-1}(1+|x-y|)^{-n-1}
\end{equation}
for all $j \in \Z_+$ and $k, \ell \in \Z^n$,
then
\begin{align*}
\|\sigma_j^{k,\ell}(X,D)f\|_{L^p}
&\le C2^{-j\delta n/p}2^{jm}(1+|k|)^{-n-1}
\|[(1+|\cdot|)^{-n-1}]*|\Lambda_{2^{-j\delta}}f|\|_{L^p} \\
&\le C2^{-j\delta n/p}2^{jm}(1+|k|)^{-n-1}
\|\Lambda_{2^{-j\delta}}f\|_{L^p} \\
&=C2^{jm}(1+|k|)^{-n-1}\|f\|_{L^p}
\end{align*}
for all $j \in \Z_+$,
$k, \ell \in \Z^n$ and $f \in \calS(\R^n)$.
This is the desired result.
We prove \eqref{(3.8)}.
Set
$\tau_j(x,\xi)
=2^{-jm}\sigma(2^{-j\delta}x,2^{j\delta}\xi)\, \psi_j(2^{j\delta}\xi)$.
Since
\[
\sigma_j^{k,\ell}(x,\xi)
=\left(\int_{\R^n}e^{ik\cdot(2^{j\delta}x-z)}\,
\Phi(2^{j\delta}x-z)\,
\sigma(2^{-j\delta}z,\xi)\, \psi_j(\xi)\,
dz\right) \varphi(2^{-j\delta}\xi-\ell),
\]
we have
\[
K_j^{k,\ell}(x,y)
=\frac{2^{jm}}{(2\pi)^n}\int_{\R^n}e^{i(x-y)\cdot\xi}\left\{
\left( \int_{\R^n}e^{ik\cdot (x-z)}\,
\Phi(x-z)\, \tau_j(z,\xi)\,
dz \right) \varphi(\xi-\ell)\right\}d\xi,
\]
where $\Phi=\calF^{-1}\varphi$.
On the other hand,
for any $\alpha,\beta \in \Z_+^n$,
\[
|\partial_{\xi}^{\alpha}\partial_x^{\beta}\tau_j(x,\xi)|
\le C_{\alpha,\beta}|\sigma|_{S_{\rho,\delta}^m, |\alpha+\beta|}
\]
for all $j \in \Z_+$.
Let $\alpha, \beta \in \Z_+^n$ be such that
$|\alpha|, |\beta| \le n+1$.
Using
\begin{align*}
&(x-y)^{\alpha}\, k^{\beta}\, K_j^{k,\ell}(x,y) \\
&=C_{\alpha,\beta}2^{jm}\int_{\R^n}
\left( \partial_{\xi}^{\alpha}e^{i(x-y)\cdot\xi}\right)
\left\{ \int_{\R^n}\left(\partial_z^{\beta}e^{ik\cdot (x-z)}\right)
\Phi(x-z)\, \tau_j(z,\xi)\,
dz \right\} \varphi(\xi-\ell)\, d\xi \\
&=C_{\alpha,\beta}2^{jm}
\sum_{\scriptstyle \alpha_1+\alpha_2=\alpha
\atop \scriptstyle \beta_1+\beta_2=\beta}
C_{\beta_1,\beta_2}^{\alpha_1,\alpha_2} \\
&\quad \times \int_{\R^n}
e^{i(x-y)\cdot\xi}
\left\{ \int_{\R^n}e^{ik\cdot (x-z)}
\partial^{\beta_1}\Phi(x-z)\,
\partial_{\xi}^{\alpha_1}\partial_z^{\beta_2}\tau_j(z,\xi)\,
dz \right\} \partial^{\alpha_2}\varphi(\xi-\ell)\, d\xi,
\end{align*}
we see that
\[
|(x-y)^{\alpha}\, k^{\beta}\, K_j^{k,\ell}(x,y)|
\le \left(C_{\alpha,\beta}\, |\sigma|_{S_{\rho,\delta}^m, |\alpha+\beta|}\,
\sup_{\tilde{\alpha} \le \alpha}
\|\partial^{\tilde{\alpha}}\varphi\|_{L^1}\, 
\sup_{\tilde{\beta} \le \beta}
\|\partial^{\tilde{\beta}}\Phi\|_{L^1}\right) 2^{jm}.
\]
This implies \eqref{(3.8)}.
The proof is complete.
\end{proof}
For $0\le \delta <1$,
we take a sufficiently large integer $j_0$
such that
\begin{equation}\label{(3.9)}
2^{j_0(1-\delta)-3} \ge \sqrt{n}.
\end{equation}
We recall that $\varphi, \psi \in \calS(\R^n)$ satisfy
$\varphi \subset [-1,1]^n$ and
$\psi \subset \{\xi: 2^{-1} \le |\xi| \le 2\}$
(see \eqref{(3.4)} and \eqref{(3.5)}).
\begin{lem}\label{3.4}
Let $0\le \delta <1$ and $\ell \in \Z^n$.
If there exists a positive integer $j(\ell) \ge j_0+1$
such that
$\mathrm{supp}\, \varphi(2^{-j(\ell)\delta}\cdot-\ell) \cap
\mathrm{supp}\, \psi(2^{-j(\ell)}\cdot) \neq \emptyset$,
then
\[
\mathrm{supp}\, \varphi(2^{-j\delta}\cdot-\ell) \cap
\mathrm{supp}\, \psi(2^{-j}\cdot)=\emptyset
\quad \text{for all} \ |j-j(\ell)| \ge j_0,
\]
where $j_0$ is defined by \eqref{(3.9)}.
\end{lem}
\begin{proof}
Since
$\mathrm{supp}\, \varphi(2^{-j\delta}\cdot-\ell)$
$\subset$
$\{|\xi-2^{j\delta}\ell| \le 2^{j\delta}\sqrt{n}\}$
and
$\mathrm{supp}\, \psi(2^{-j}\cdot)$
$\subset$
$\{2^{j-1}\le |\xi| \le 2^{j+1}\}$,
our assumption
$\mathrm{supp}\, \varphi(2^{-j(\ell)\delta}\cdot-\ell)$
$\cap$
$\mathrm{supp}\, \psi(2^{-j(\ell)}\cdot)$
$\neq$
$\emptyset$
gives
\begin{equation}\label{(3.10)}
2^{j(\ell)-1}-2^{j(\ell)\delta}\sqrt{n}
\le 2^{j(\ell)\delta}|\ell|
\le 2^{j(\ell)+1}+2^{j(\ell)\delta}\sqrt{n}.
\end{equation}
Let $|j-j(\ell)| \ge j_0$.
We consider the case $j-j(\ell) \le -j_0$.
By \eqref{(3.9)},
\eqref{(3.10)} and $j(\ell) \ge j_0+1$,
we see that
\begin{align*}
&2^{j\delta}|\ell|-2^{j\delta}\sqrt{n}
=2^{j\delta-j(\ell)\delta+j(\ell)\delta}|\ell|-2^{j\delta}\sqrt{n} \\
&\ge 2^{j\delta-j(\ell)\delta}(2^{j(\ell)-1}-2^{j(\ell)\delta}\sqrt{n})
-2^{j\delta}\sqrt{n}
=2^{j\delta+j(\ell)(1-\delta)-1}-2^{j\delta+1}\sqrt{n} \\
&\ge 2^{j\delta+j(\ell)(1-\delta)-1}-2^{j\delta+j_0(1-\delta)-2}
>2^{j\delta+j(\ell)(1-\delta)-1}-2^{j\delta+j(\ell)(1-\delta)-2} \\
&=2^{j\delta+j(\ell)(1-\delta)-2}
=2^{j\delta+j(\ell)(1-\delta)-2-(j+1)+(j+1)}
=2^{\{(j(\ell)-j)(1-\delta)-3\}+j+1} \\
&\ge 2^{(j_0(1-\delta)-3)+j+1}
\ge 2^{j+1}\sqrt{n} \ge 2^{j+1}.
\end{align*}
Hence,
\[
\mathrm{supp}\, \varphi(2^{-j\delta}\cdot-\ell)
\subset \{|\xi| \ge 2^{j\delta}|\ell|-2^{j\delta}\sqrt{n}\}
\subset \{|\xi|>2^{j+1}\}
\]
for all $j-j(\ell) \le -j_0$.
This implies
$\mathrm{supp}\, \varphi(2^{-j\delta}\cdot-\ell) \cap
\mathrm{supp}\, \psi(2^{-j}\cdot)=\emptyset$
for all $j-j(\ell) \le -j_0$.

In the same way, we can prove
\[
\mathrm{supp}\, \varphi(2^{-j\delta}\cdot-\ell)
\subset \{|\xi| \le 2^{j\delta}|\ell|+2^{j\delta}\sqrt{n}\}
\subset \{|\xi|<2^{j-1}\}
\]
for all $j-j(\ell) \ge j_0$.
This implies
$\mathrm{supp}\, \varphi(2^{-j\delta}\cdot-\ell) \cap
\mathrm{supp}\, \psi(2^{-j}\cdot)=\emptyset$
for all $j-j(\ell) \ge j_0$.
\end{proof}
Let $\eta \in \calS(\R^n)$ be such that
$\mathrm{supp}\, \eta$ is compact and
$|\sum_{\nu \in \Z^n}\eta(\xi-\nu)| \ge C>0$
for all $\xi \in \R^n$.
It is well known that
\begin{equation}\label{(3.11)}
\|f\|_{M^{p,q}} \sim
\left(\sum_{\nu \in \Z^n}\|\eta(D-\nu)f\|_{L^p}^q\right)^{1/q},
\end{equation}
where $\eta(D-\nu)f=\calF^{-1}[\eta(\cdot-\nu)\widehat{f}]$
(see, for example, \cite{Triebel}).
\begin{prop}\label{3.5}
Let $1 \le p \le \infty$,
$m \in \R$, $0 \le \delta \le \rho \le 1$,
$\delta<1$ and $\sigma \in S_{\rho,\delta}^m$.
If $m \le -(\mu_1(p,\infty)-\mu_2(p,\infty))\delta n$,
then there exists a constant $C>0$ such that
\[
\|\sigma(X,D)f\|_{M^{p,\infty}} \le C\|f\|_{M^{p,\infty}}
\]
for all $f \in \calS(\R^n)$.
\end{prop}
\begin{proof}
In view of \eqref{(3.11)}, we estimate
$\sup_{\nu \in \Z^n}\|\varphi(D-\nu)(\sigma(X,D)f)\|_{L^p}$,
where $\varphi$ is as \eqref{(3.5)}
and $f \in \calS(\R^n)$.
By the decomposition \eqref{(3.6)},
we have
\begin{align*}
\sigma(x,\xi)
&=\sum_{k \in \Z^n}\sum_{\ell \in \Z^n}\sum_{j=0}^{\infty}
\sigma_j^{k,\ell}(x,\xi) \\
&=\sum_{k \in \Z^n}\sum_{\ell \in \Z^n}
\sum_{j=0}^{j_0}
\sigma_j^{k,\ell}(x,\xi)
+\sum_{k \in \Z^n}\sum_{\ell \in \Z^n}
\sum_{j=j_0+1}^{\infty}
\sigma_j^{k,\ell}(x,\xi),
\end{align*}
where $\sigma_j^{k,\ell}(x,\xi)$ is defined by \eqref{(3.7)},
and $j_0$ is defined by \eqref{(3.9)}.
We only consider the second sum since
the estimate for the first sum can be carried out in a similar way.
By Lemma \ref{3.4},
we see that
\[
\sum_{k \in \Z^n}\sum_{\ell \in \Z^n}
\sum_{j=j_0+1}^{\infty}
\sigma_j^{k,\ell}(x,\xi)
=\sum_{k \in \Z^n}\sum_{\ell \in \Z^n}
\sum_{\scriptstyle j\ge j_0+1 \scriptstyle \atop |j-j(\ell)|<j_0}
\sigma_j^{k,\ell}(x,\xi).
\]
Since
$\varphi(2^{-j\delta}\cdot-k)
\subset 2^{j\delta}k+[-2^{j\delta},2^{j\delta}]^n$
and
\begin{align*}
&\calF[\sigma_j^{k,\ell}(X,D)f](y) \\
&=\frac{1}{(2\pi)^n}
\int_{\R^n}
\calF_{x \to y}[e^{ix\cdot\xi}
\varphi(2^{-j\delta}D_x-k)\,\sigma(x,\xi)]\,
\varphi(2^{-j\delta}\xi-\ell)\, \psi(2^{-j}\xi)\, \widehat{f}(\xi)\, d\xi \\
&=\frac{1}{(2\pi)^n}
\int_{\R^n}
\varphi(2^{-j\delta}(y-\xi)-k)\, \calF_1\sigma(y-\xi,\xi)\, 
\varphi(2^{-j\delta}\xi-\ell)\, \psi(2^{-j}\xi)\, \widehat{f}(\xi)\, d\xi,
\end{align*}
we see that
$\calF[\sigma_j^{k,\ell}(X,D)f]
\subset 2^{j\delta}(k+\ell)+[-2^{j\delta+1},2^{j\delta+1}]^n$.
This gives
\begin{align*}
&\varphi(D-\nu)\left(\sum_{k \in \Z^n}\sum_{\ell \in \Z^n}
\sum_{\scriptstyle j\ge j_0+1 \scriptstyle \atop |j-j(\ell)|<j_0}
\sigma_j^{k,\ell}(X,D)f \right) \\
&=\sum_{k \in \Z^n}
\sum_{\scriptstyle \ell \in \Z^n, i=1,\dots, n
\scriptstyle \atop |2^{j\delta}(k_i+\ell_i)-\nu_i| \le 2^{j\delta+2}}
\sum_{\scriptstyle j\ge j_0+1 \scriptstyle \atop |j-j(\ell)|<j_0}
\varphi(D-\nu)\, \sigma_j^{k,\ell}(X,D)f.
\end{align*}
Let $\eta \in \calS(\R^n)$ be such that
$\mathrm{supp}\, \eta$ is compact,
$\eta=1$ on $\mathrm{supp}\, \varphi$ and
$|\sum_{\nu \in \Z^n}\eta(\xi-\nu)| \ge C>0$
for all $\xi \in \R^n$.
Note that
$\sigma_j^{k,\ell}(x,\xi)
=\sigma_j^{k,\ell}(x,\xi)\, \eta(2^{-j\delta}\xi-\ell)$.
By Lemma \ref{3.3},
we have
\begin{align*}
&\left\|\varphi(D-\nu)\left(\sum_{k \in \Z^n}\sum_{\ell \in \Z^n}
\sum_{j_0+1}^{\infty}
\sigma_j^{k,\ell}(X,D)f \right)\right\|_{L^p} \\
&\le \sum_{k \in \Z^n}
\sum_{\scriptstyle \ell \in \Z^n, i=1,\dots, n
\scriptstyle \atop |2^{j\delta}(k_i+\ell_i)-\nu_i| \le 2^{j\delta+2}}
\sum_{\scriptstyle j\ge j_0+1 \scriptstyle \atop |j-j(\ell)|<j_0}
\|\varphi(D-\nu)\, \sigma_j^{k,\ell}(X,D)f\|_{L^p} \\
&\le C\sum_{k \in \Z^n}
\sum_{\scriptstyle \ell \in \Z^n, i=1,\dots, n
\scriptstyle \atop |2^{j\delta}(k_i+\ell_i)-\nu_i| \le 2^{j\delta+2}}
\sum_{\scriptstyle j\ge j_0+1 \scriptstyle \atop |j-j(\ell)|<j_0}
\|\sigma_j^{k,\ell}(X,D)[\eta(2^{-j\delta}D-\ell)f]\|_{L^p} \\
&\le C\sum_{k \in \Z^n}
\sum_{\scriptstyle \ell \in \Z^n, i=1,\dots, n
\scriptstyle \atop |2^{j\delta}(k_i+\ell_i)-\nu_i| \le 2^{j\delta+2}}
\sum_{\scriptstyle j\ge j_0+1 \scriptstyle \atop |j-j(\ell)|<j_0}
2^{jm}(1+|k|)^{-n-1}\|\eta(2^{-j\delta}D-\ell)f\|_{L^p}.
\end{align*}
Since
$\eta(2^{-j\delta}D-\ell)f
=\Lambda_{2^{j\delta}}\, \eta(D-\ell)\, \Lambda_{2^{-j\delta}}f$,
$\mu_1(p,\infty)=-1/p$ and
$m+(\mu_1(p,\infty)-\mu_2(p,\infty))\delta n \le 0$,
by \eqref{(3.11)} and Proposition \ref{2.1},
we see that
\begin{align*}
&\sum_{k \in \Z^n}
\sum_{\scriptstyle \ell \in \Z^n, i=1,\dots, n
\scriptstyle \atop |2^{j\delta}(k_i+\ell_i)-\nu_i| \le 2^{j\delta+2}}
\sum_{\scriptstyle j\ge j_0+1 \scriptstyle \atop |j-j(\ell)|<j_0}
2^{jm}(1+|k|)^{-n-1}\|\eta(2^{-j\delta}D-\ell)f\|_{L^p} \\
&=\sum_{k \in \Z^n}
\sum_{\scriptstyle \ell \in \Z^n, i=1,\dots, n
\scriptstyle \atop |\ell_i-(2^{-j\delta}\nu_i-k_i)| \le 4}
\sum_{\scriptstyle j\ge j_0+1 \scriptstyle \atop |j-j(\ell)|<j_0}
2^{j(m-\delta n/p)}(1+|k|)^{-n-1}
\|\eta(D-\ell)(\Lambda_{2^{-j\delta}}f)\|_{L^p} \\
&\le C \sum_{k \in \Z^n}
\sum_{\scriptstyle \ell \in \Z^n, i=1,\dots, n
\scriptstyle \atop |\ell_i-(2^{-j\delta}\nu_i-k_i)| \le 4}
\sum_{\scriptstyle j\ge j_0+1 \scriptstyle \atop |j-j(\ell)|<j_0}
2^{j(m-\delta n/p)}(1+|k|)^{-n-1}
\|\Lambda_{2^{-j\delta}}f\|_{M^{p,\infty}} \\
&\le C\sum_{k \in \Z^n}
\sum_{\scriptstyle \ell \in \Z^n, i=1,\dots, n
\scriptstyle \atop |\ell_i-(2^{-j\delta}\nu_i-k_i)| \le 4}
\sum_{\scriptstyle j\ge j_0+1 \scriptstyle \atop |j-j(\ell)|<j_0}
2^{j(m+(\mu_1(p,\infty)-\mu_2(p,\infty))\delta n)}(1+|k|)^{-n-1}
\|f\|_{M^{p,\infty}} \\
&\le C\left\{\sum_{k \in \Z^n}(1+|k|)^{-n-1}
\sum_{\scriptstyle \ell \in \Z^n, i=1,\dots, n
\scriptstyle \atop |\ell_i-(2^{-j\delta}\nu_i-k_i)| \le 4}
\left(\sum_{\scriptstyle j\ge j_0+1 \scriptstyle \atop |j-j(\ell)|<j_0}1
\right) \right\}\|f\|_{M^{p,\infty}} \\
&\le C(2j_0-1)
\left\{\sum_{k \in \Z^n}(1+|k|)^{-n-1}
\left(\sum_{\scriptstyle \ell \in \Z^n, i=1,\dots, n
\scriptstyle \atop |\ell_i-(2^{-j\delta}\nu_i-k_i)| \le 4}1\right)\right\}
\|f\|_{M^{p,\infty}} \\
&\le C(2j_0-1)9^n
\left(\sum_{k \in \Z^n}(1+|k|)^{-n-1}\right)
\|f\|_{M^{p,\infty}}
=C_{n, \delta}\|f\|_{M^{p,\infty}}.
\end{align*}
This implies
\begin{align*}
&\left\|\sum_{k \in \Z^n}\sum_{\ell \in \Z^n}
\sum_{j_0+1}^{\infty}
\sigma_j^{k,\ell}(X,D)f\right\|_{M^{p,\infty}} \\
&\le C\sup_{\nu \in \Z^n}
\left\|\varphi(D-\nu)\left(\sum_{k \in \Z^n}\sum_{\ell \in \Z^n}
\sum_{j_0+1}^{\infty}
\sigma_j^{k,\ell}(X,D)f \right)\right\|_{L^p}
\le C_{n,\delta}\|f\|_{M^{p,\infty}}
\end{align*}
for all $f \in \calS(\R^n)$.
This is the desired result.
\end{proof}
In order to prove Theorem \ref{3.1},
we will use the following interpolation technique
(see Toft \cite[Remark 3.2]{Toft}).
Let $\calM^{p,q}(\R^n)$ be the completion of $\calS(\R^n)$
under the norm $\|\cdot\|_{M^{p,q}}$.
Then the following are true:
\begin{enumerate}
\item
If $1 \le p,q<\infty$, then $\calM^{p,q}(\R^n)=M^{p,q}(\R^n)$,
\item
If $1 \le p_1,p_2,q_1,q_2 \le \infty$,
$0 < \theta <1$,
$1/p=\theta/p_1+(1-\theta)/p_2$
and $1/q=\theta/q_1+(1-\theta)/q_2$,
then
$(\calM^{p_1,q_1}(\R^n), \calM^{p_2,q_2}(\R^n))_{\theta}
=\calM^{p,q}(\R^n)$.
\end{enumerate}

We are now ready to prove Theorem \ref{3.1}.

\medskip
\noindent
{\it Proof of Theorem \ref{3.1}.}
Let $1<p,q<\infty$, $0\le \delta \le \rho \le 1$
and $\delta<1$.
We first consider the case $(1/p,1/q) \in J_1$ with $1/p \ge 1/2$.
Then $\mu_1(p,q)-\mu_2(p,q)=1/p-1/q$.
We take $1 \le r \le 2$ and $0 < \theta < 1$
such that $1/p=\theta/2+(1-\theta)/r$
and $1/q=\theta/2+(1-\theta)/\infty$.
Since $\mu_1(r,\infty)-\mu_2(r,\infty)=1/r$,
by Proposition \ref{3.5},
we have
\begin{equation}\label{(3.12)}
\mathrm{Op}(S_{\rho,\delta}^{-\delta n/r})
\subset \calL(\calM^{r,\infty}(\R^n)).
\end{equation}
On the other hand,
it is well known that
$\mathrm{Op}(S_{\rho,\delta}^0) \subset \calL(L^2(\R^n))$
(\cite[Chapter 7, Theorem 2]{Stein}).
Hence,
\begin{equation}\label{(3.13)}
\mathrm{Op}(S_{\rho,\delta}^0)
\subset \calL(\calM^{2,2}(\R^n)).
\end{equation}
By interpolation,
\eqref{(3.12)} and \eqref{(3.13)} give
$\mathrm{Op}(S_{\rho,\delta}^{0\cdot\theta+(-\delta n/r)(1-\theta)})
\subset \calL(\calM^{p,q}(\R^n))$.
Since $(1-\theta)/r=1/p-1/q$ and $\calM^{p,q}(\R^n)=M^{p,q}(\R^n)$,
we obtain
\begin{equation}\label{(3.14)}
\mathrm{Op}(S_{\rho,\delta}^{-(1/p-1/q)\delta n})
\subset \calL(M^{p,q}(\R^n))
\quad \text{if $(1/p,1/q) \in J_1$ with $1/p \ge 1/2$}.
\end{equation}
In the same way,
we can prove
\begin{equation}\label{(3.15)}
\mathrm{Op}(S_{\rho,\delta}^{-(-1/p-1/q+1)\delta n})
\subset \calL(M^{p,q}(\R^n))
\quad \text{if $(1/p,1/q) \in J_3$ with $1/p \le 1/2$}.
\end{equation}
Note that $\mu_1(p,q)-\mu_2(p,q)=-1/p-1/q+1$
if $(1/p,1/q) \in J_3$ with $1/p \le 1/2$.

We next consider the case $(1/p,1/q) \in J_1$ with $1/p \le 1/2$.
Then $\mu_1(p,q)-\mu_2(p,q)=-1/p+1/q$.
Let $\sigma \in S_{\rho,\delta}^{-(-1/p+1/q)\delta n}$.
Then $\sigma^*(X,D)$ satisfies \eqref{(3.2)},
where $\sigma^* \in S_{\rho,\delta}^{-(-1/p+1/q)\delta n}$
is defined by \eqref{(3.1)}.
Let $p',q'$ be the conjugate exponents of $p,q$,
respectively.
Since $-(-1/p+1/q)\delta n=-(1/p'-1/q')\delta n$
and $(1/p',1/q') \in J_1$ with $1/p' \ge 1/2$,
by \eqref{(3.14)},
we see that $\sigma^*(X,D)$ is bounded on $M^{p',q'}(\R^n)$.
Then,
by duality and \eqref{(3.2)},
we obtain that $\sigma(X,D)$ is bounded on $M^{p,q}(\R^n)$,
that is,
\begin{equation}\label{(3.16)}
\mathrm{Op}(S_{\rho,\delta}^{-(-1/p+1/q)\delta n})
\subset \calL(M^{p,q}(\R^n))
\quad \text{if $(1/p,1/q) \in J_1$ with $1/p \le 1/2$}.
\end{equation}
In the same way,
using duality and \eqref{(3.15)},
we can prove
\begin{equation}\label{(3.17)}
\mathrm{Op}(S_{\rho,\delta}^{-(1/p+1/q-1)\delta n})
\subset \calL(M^{p,q}(\R^n))
\quad \text{if $(1/p,1/q) \in J_3$ with $1/p \ge 1/2$}.
\end{equation}
Note that $\mu_1(p,q)-\mu_2(p,q)=1/p+1/q-1$
if $(1/p,1/q) \in J_3$ with $1/p \ge 1/2$.

Finally, we consider the case $(1/p,1/q) \in J_2$.
Let $1/p \le 1/2$.
Then $\mu_1(p,q)-\mu_2(p,q)=-2/p+1$.
Since $(1/p,1/p') \in J_1$ with $1/p \le 1/2$
and $(1/p,1/p) \in J_3$ with $1/p \le 1/2$,
\eqref{(3.16)} and \eqref{(3.15)} imply
$\mathrm{Op}(S_{\rho,\delta}^{-(-1/p+1/p')\delta n})
\subset \calL(M^{p,p'}(\R^n))$
and $\mathrm{Op}(S_{\rho,\delta}^{-(-2/p+1)\delta n})
\subset \calL(M^{p,p}(\R^n))$.
Then,
by interpolation,
we have
\[\mathrm{Op}
(S_{\rho,\delta}^{-(-1/p+1/p')\delta n \theta-(-2/p+1)\delta n(1-\theta)})
\subset \calL(M^{p,q}(\R^n)),
\]
where $1/q=\theta/p'+(1-\theta)/p$,
that is,
\begin{equation}\label{(3.18)}
\mathrm{Op}(S_{\rho,\delta}^{-(-2/p+1)\delta n})
\subset \calL(M^{p,q}(\R^n))
\quad \text{if $(1/p,1/q) \in J_2$ with $1/p \le 1/2$}.
\end{equation}
In the same way
(or duality with \eqref{(3.18)}),
using interpolation, \eqref{(3.14)} and \eqref{(3.17)},
we can prove
\[
\mathrm{Op}(S_{\rho,\delta}^{-(2/p-1)\delta n})
\subset \calL(M^{p,q}(\R^n))
\quad \text{if $(1/p,1/q) \in J_2$ with $1/p \ge 1/2$}.
\]
Note that $\mu_1(p,q)-\mu_2(p,q)=2/p-1$
if $(1/p,1/q) \in J_2$ with $1/p \ge 1/2$.
The proof is complete.

\medskip
In order to prove the $\lq\lq$only if'' part of Theorem \ref{1.1},
we introduce a special symbol.
Let $\varphi,\eta \in \calS(\R^n)$ be real-valued functions
satisfying
\begin{align*}
\varphi:
\quad &{\rm supp}\, \varphi \subset \{\xi: |\xi| \le 1/8\},
\quad \int_{\R^n}\varphi(\xi)\, d\xi =1, \\
\eta:
\quad &{\rm supp}\, \eta \subset
\{\xi: 2^{-1/2}\le |\xi| \le 2^{1/2}\},
\quad \eta=1 \ \text{on} \
\{\xi: 2^{-1/4}\le |\xi| \le 2^{1/4}\}.
\end{align*}
Moreover, we assume that $\varphi$ is radial.
This assumption implies that $\Phi$ is also real-valued,
where $\Phi=\calF^{-1}\varphi$.
Then we define
\begin{equation}\label{(3.19)}
\sigma(x,\xi)
=\sum_{j=j_0}^{\infty}2^{jm}
\left(\sum_{0<|k|\le 2^{j\delta/2}}
e^{-ik\cdot(2^{j\delta/2}x-k)}\,
\Phi(2^{j\delta/2}x-k) \right)
\eta(2^{-j}\xi),
\end{equation}
where $0<\delta<1$ and $j_0\in\Z_+$ is chosen to satisfy
\begin{equation}\label{(3.20)}
1+2^{j_0(\delta-1)+1} \le 2^{1/4}, \quad
1-2^{j_0(\delta-1)+1} \ge 2^{-1/4}, \quad
2^{-j_0 \delta/2} \sqrt{n} \le 2^{-3}.
\end{equation}
The symbol $\sigma^*$ is constructed from
$\sigma$ using the oscillatory integral \eqref{(3.1)}.
The following is the $\lq\lq$only if'' part of Theorem \ref{1.1}.
\begin{thm}[{\cite[Theorem 1.1]{Sugimoto-Tomita-2}}]\label{3.6}
Let $1<p,q<\infty$, $0\le\delta<1$ and $m>-|1/q-1/2|\delta n$.
Then the symbols $\sigma$ and $\sigma^*$
defined by \eqref{(3.19)}
belong to the class $S_{1,\delta}^m$.
Moreover,
if $q\geq2$ {\rm (}$q\leq2$ resp.{\rm )}, then the corresponding operator
$\sigma(X,D)$ {\rm (}$\sigma^*(X,D)$ resp.{\rm )} is
not bounded on $M^{p,q}(\R^n)$.
\end{thm}
Precisely speaking,
\cite{Sugimoto-Tomita-2} treated only the case $0<\delta<1$,
but the case $\delta=0$ can be also included by a simple argument.
If $m>0$ and
$\sigma(\xi)=(1+|\xi|^2)^{m/2}$,
then $\sigma(D)$ is not bounded on $L^2(\R^n)$.
Indeed,
\begin{align*}
\|\sigma(D)(M_kf)\|_{L^2}
&=(2\pi)^{-n/2}\|\widehat{\sigma(D)(M_kf)}\|_{L^2} \\
&=(2\pi)^{-n/2}\|(1+|\xi|^2)^{m/2}\, \widehat{f}(\cdot-k)\|_{L^2}
\ge C|k|^m
\end{align*}
for all $k \in \Z^n$,
where $f \in \calS(\R^n)\setminus \{0\}$ such that
${\rm supp}\, \widehat{f} \subset \{|\xi| \le 1/2\}$.
On the other hand,
$\|M_kf\|_{L^2}=\|f\|_{L^2}$
for all $k \in \Z^n$.
Hence, $\sigma(D)$ is not bounded on $L^2(\R^n)$.
Moreover, $\sigma(D)$ is not bounded on $M^{p,q}(\R^n)$
for any $1<p,q<\infty$.
Indeed, since $\sigma(D)^*=\sigma(D)$,
if $\sigma(D)$ is bounded on $M^{p,q}(\R^n)$,
then $\sigma(D)$ is also bounded on $M^{p',q'}(\R^n)$.
Then, by interpolation,
$\sigma(D)$ is bounded on $M^{2,2}(\R^n)$.
However, since $M^{2,2}(\R^n)=L^2(\R^n)$,
this is a contradiction.
Hence, $\sigma(D)$ is not bounded on $M^{p,q}(\R^n)$
for any $1<p,q<\infty$.
Note that $\sigma(\xi)=(1+|\xi|^2)^{m/2}$
belongs to $S_{1,0}^m$.

\section{The unboundedness of Calder\'on-Zygmund operators}\label{section4}
We recall the theory of Calder\'on-Zygmund operators
based on \cite{F-T-W}
(see also \cite{Torres}).
Let $\calD(\R^n)$ be the space
of all infinitely differentiable functions with compact support.
For $\ell \in \Z_+$,
we denote by $\calD_{\ell}(\R^n)$
the space of all $g \in \calD(\R^n)$ such that
$\int_{\R^n}x^{\beta}g(x)\, dx=0$
for all $|\beta| \le \ell$.
Let $T: \calD(\R^n) \to \calD'(\R^n)$ be
a continuous linear operator,
where $\calD'(\R^n)$ is the dual space of $\calD(\R^n)$.
We denote by $K$ the distributional kernel of $T$,
that is,
\[
\langle Tf, g \rangle
=\langle K, g\otimes f \rangle
\qquad \text{for all}\ f,g \in \calD(\R^n).
\]
We say that $T$ is a generalized Calder\'on-Zygmund operator
of smoothness $\ell+\epsilon$
(we write $T \in {\rm CZO}(\ell+\epsilon)$),
where $\ell \in \Z_+$ and $0<\epsilon<1$,
if the restriction of $K$
to $\{(x,y)\in \R^n \times \R^n: x\neq y\}$
is a function
with continuous partial derivatives
in the variable $x$ up to order $\ell$
which satisfy
\begin{align*}
&|\partial_x^{\alpha}K(x,y)| \le
C|x-y|^{-n-|\alpha|}
\qquad \text{for} \ |\alpha| \le \ell, \\
&|\partial_x^{\alpha}K(x,y)-\partial_x^{\alpha}K(x',y)|
\le C\frac{|x-x'|^{\epsilon}}{|x-y|^{n+\ell+\epsilon}}
\end{align*}
for $|\alpha|=\ell$ and $|x-y| >2|x-x'|$.
Let $\phi \in \calD(\R^n)$ be such that
$\phi=1$ on $\{x:|x| \le 1\}$.
If $T \in {\rm CZO}(\ell+\epsilon)$
and $f \in C^{\infty}(\R^n)$ satisfies $f(x)=O(|x|^\ell)$
as $|x| \to \infty$,
then the limit
$\lim_{j \to \infty}\langle T(\phi(\cdot/j)f),g \rangle$
exists for all $g \in \calD_\ell(\R^n)$
(\cite[Lemma 1.19]{F-T-W}, \cite[Lemma 2.2.12]{Torres}).
In particular,
$T(x^{\beta})$ can be defined as an element of $\calD_{\ell}'(\R^n)$,
where $|\beta|\le \ell$
and $\calD_{\ell}'(\R^n)$ is the dual space of $\calD_{\ell}(\R^n)$.
Note that this limit
$\lim_{j \to \infty}\langle T(\phi(\cdot/j)f),g \rangle$
is independent of the choice of $\phi$
(\cite[p.49]{F-T-W}).
It is known that,
if $T \in {\rm CZO}(\ell+\epsilon)$
is a translation invariant operator,
then $T(x^{\beta})=0$ as an element of $\calD_{k}'(\R^n)$
for all $|\beta|=k\le \ell$
(\cite[Proposition 2.2.17]{Torres}).
A continuous linear operator $T$ is said to satisfy
the weak boundedness property,
if for each bounded subset $\calB$ of $\calD(\R^n)$
there exists a constant $C=C(\calB)>0$ such that
\[
|\langle Tf^{x,R},g^{x,R} \rangle |
\le CR^n
\qquad \text{for all}\ f,g \in \calB,
\ x \in \R^n
\ \text{and} \ R>0,
\]
where $f^{x,R}(y)=f((y-x)/R)$.
The transpose of $T$,
$T^*:\calD(\R^n) \to \calD'(\R^n)$
is defined by
\[
\langle T^*f,g \rangle
= \langle Tg, f \rangle
\qquad \text{for} \ f,g \in \calD(\R^n).
\] 
Note that the kernel of $T^*$ is given by
$K^*(x,y)=K(y,x)$.
David and Journ\'e \cite{D-J} proved that,
when $T,T^* \in {\rm CZO}(\epsilon)$,
$T$ is bounded on $L^2(\R^n)$
if and only if
$T(1), T^*(1) \in BMO(\R^n)$
and $T$ satisfies the weak boundedness property.

Our main result on the unboundedness of Calder\'on-Zygmund operators
is the following:
\begin{thm}\label{4.1}
Let $1<p,q<\infty$.
If $q \neq 2$,
then there exists an operator $T: \calS(\R^n) \to \calS'(\R^n)$
such that
$T, T^* \in {\rm CZO}(\ell+\epsilon)$
for all $\ell \in \Z_+$,
where $0<\epsilon<1$,
$T(P)=T^*(P)=0$
for all polynomials $P$,
$T$ satisfies the weak boundedness property,
but $T$ is not bounded on $M^{p,q}(\R^n)$.
\end{thm}
By David-Journ\'e's $T(1)$ theorem mentioned above,
as a corollary of Theorem \ref{4.1},
we have Theorem \ref{1.2}.
In the rest of this section,
we prove Theorem \ref{4.1}.
\begin{lem}\label{4.2}
Let $\sigma$ be defined by {\rm \eqref{(3.19)}}
with $0<\delta<1$ and $j_0$ satisfying \eqref{(3.20)} and
$2^{j_0(\delta-1)+2}$
$<$
$2^{-1/2}$.
Then $\sigma, \sigma^* \in S_{1,\delta}^m$ and
$(\partial_{\xi}^{\alpha}\sigma)(x,0)
=(\partial_{\xi}^{\alpha}\sigma^*)(x,0)=0$
for all $\alpha \in \Z_+^n$,
where $\sigma^*$ is defined by {\rm \eqref{(3.1)}}.
\end{lem}
\begin{proof}
By Theorem \ref{3.6},
we have $\sigma, \sigma^* \in S_{1,\delta}^m$.
Since
${\rm supp}\, \sigma \subset \{(x,\xi): |\xi| \ge 2^{j_0-1/2}\}$,
we obtain
$(\partial_{\xi}^{\alpha}\sigma)(x,0)=0$
for all $\alpha \in \Z_+^n$.

Let $\chi_1,\chi_2 \in \calS(\R^n)$
be such that $\chi_1(0)=\chi_2(0)=1$
and
${\rm supp}\, \widehat{\chi_1}, {\rm supp}\, \chi_2
\subset \{|\zeta| \le 1\}$.
Since
${\rm supp}\, \eta (2^{-j}\cdot)
\subset \{2^{j-1/2}\le |\zeta| \le 2^{j+1/2}\}$
and
${\rm supp}\, \chi_2(\epsilon \cdot)
\subset \{|\zeta| \le 1/\epsilon\}$,
for each $0<\epsilon<1$,
there exists $N_{\epsilon}$ such that
${\rm supp}\, \eta(2^{-j}\cdot) \cap
{\rm supp}\, \chi_2(\epsilon \cdot) =\emptyset$
for all $j \ge N_{\epsilon}$.
Hence, by \eqref{(3.3)}, we see that
\begin{align*}
&(\partial_{\xi}^{\alpha}\sigma^*)(x,0)
=\textrm{Os-}
\frac{1}{(2\pi)^n}\int_{\R^n}\int_{\R^n}
e^{-iy\cdot\zeta}\,
\overline{(\partial_{\xi}^{\alpha}\sigma)(x+y,\zeta)}
\, dy \, d\zeta \\
&=\lim_{\epsilon \to 0}
\frac{1}{(2\pi)^n}
\int_{\R^n}\int_{\R^n}
e^{-iy\cdot\zeta}\,
\chi_1(\epsilon y)\, \chi_2(\epsilon \zeta)
\Bigg\{ \sum_{j=j_0}^{\infty}2^{j(m-|\alpha|)} \\
&\quad \times 
\Bigg(\sum_{0<|k|\le 2^{j\delta/2}}
e^{ik\cdot(2^{j\delta/2}(x+y)-k)} \,
\Phi(2^{j\delta/2}(x+y)-k) \Bigg)
(\partial^{\alpha}\eta)(2^{-j}\zeta) \Bigg\}
dy\, d\zeta \\
&=\lim_{\epsilon \to 0}
\sum_{j=j_0}^{\infty}2^{j(m-|\alpha|)}\sum_{0<|k|\le 2^{j\delta/2}}
\frac{1}{(2\pi)^n}\int_{\R^n}\int_{\R^n}
e^{-iy\cdot\zeta}\,
\chi_1(\epsilon y)\, \chi_2(\epsilon \zeta) \\
& \qquad \qquad \times
e^{ik\cdot(2^{j\delta/2}(x+y)-k)} \,
\Phi(2^{j\delta/2}(x+y)-k)\,
(\partial^{\alpha}\eta)(2^{-j}\zeta)\,
dy\, d\zeta \\
&=\lim_{\epsilon \to 0}
\sum_{j=j_0}^{\infty}2^{j(m-|\alpha|)}\sum_{0<|k|\le 2^{j\delta/2}}
\frac{1}{(2\pi)^n}\int_{\R^n}
\chi_2(\epsilon \zeta)\, (\partial^{\alpha}\eta)(2^{-j}\zeta) \\
& \qquad \qquad \times
\left( \int_{\R^n}e^{-iy\cdot\zeta}\,
\chi_1(\epsilon y)\,
e^{ik\cdot(2^{j\delta/2}(x+y)-k)}\,
\Phi(2^{j\delta/2}(x+y)-k)\,
dy\right) d\zeta \\
&=\lim_{\epsilon \to 0}
\sum_{j=j_0}^{\infty}2^{j(m-|\alpha|)}\sum_{0<|k|\le 2^{j\delta/2}}
\frac{1}{(2\pi)^n}\int_{\R^n}
\chi_2(\epsilon \zeta)\, (\partial^{\alpha}\eta)(2^{-j}\zeta)\,
\eta_{j,k}^{\epsilon,\delta}(x,\zeta)\, d\zeta.
\end{align*}
By Plancherel's theorem,
we have
\begin{align*}
&\eta_{j,k}^{\epsilon,\delta}(x,\zeta)
=\int_{\R^n}e^{-iy\cdot\zeta}\,
\chi_1(\epsilon y)\,
\overline{e^{-ik\cdot(2^{j\delta/2}(x+y)-k)}\,
\Phi(2^{j\delta/2}(x+y)-k)}\, dy \\
&=\frac{1}{(2\pi)^n}
\int_{\R^n}
\epsilon^{-n}\widehat{\chi_1}((\zeta+\zeta')/\epsilon)\,
\overline{ 2^{-j\delta n/2}\,
e^{ix\cdot \zeta'}\, e^{-ik\cdot(2^{-j\delta/2}\zeta')}\,
\varphi(2^{-j\delta/2}\zeta'+k)} \, d\zeta'.
\end{align*}
Let $j \ge j_0$, $0<|k| \le 2^{j\delta/2}$,
$0<\epsilon<1$ and $x \in \R^n$.
Since
${\rm supp}\, \widehat{\chi_1}(\cdot/\epsilon)
\subset \{|\zeta'| \le 1\}$
and
${\rm supp}\, \varphi(2^{-j\delta/2}\cdot+k)
\subset \{|\zeta'| \le 2^{j\delta+1}\}$,
we see that
${\rm supp}\, \eta_{j,k}^{\epsilon,\delta}(x, \cdot)
\subset \{|\zeta| \le 2^{j\delta+1}+1\}$.
Therefore,
our assumption $2^{j_0(\delta-1)+2}<2^{-1/2}$
gives
${\rm supp}\, \eta_{j,k}^{\epsilon,\delta}(x,\cdot)
\subset \{|\xi|<2^{j-1/2}\}$.
This implies
${\rm supp}\, (\partial^{\alpha}\eta)(2^{-j}\cdot) \cap
{\rm supp}\, \eta_{j,k}^{\epsilon,\delta}(x,\cdot)=\emptyset$,
that is,
$(\partial_{\xi}^{\alpha}\sigma^*)(x,0)=0$.
The proof is complete.
\end{proof}
\begin{lem}[{\cite[Chapter 7, Section 1, Proposition 1]{Stein}},
{\cite[Lemma 5.1.6]{Torres}}]\label{4.3}
Let $\sigma \in S_{1,1}^{0}$.
Then, $K$, the distributional kernel of $\sigma(X,D)$
as an operator from $\calS(\R^n)$ to $\calS'(\R^n)$,
is a $C^{\infty}$-function
on $\{(x,y)\in \R^n \times \R^n: x \neq y\}$
satisfying
\[
|\partial_x^{\alpha}\partial_y^{\beta}K(x,y)|
\le C_{\alpha,\beta}|x-y|^{-n-|\alpha|-|\beta|}
\]
for all $\alpha,\beta \in \Z_+^n$.
\end{lem}
Lemma \ref{4.3} says that,
if $\sigma \in S_{1,1}^{0}$,
then $\sigma(X,D), \sigma(X,D)^* \in {\rm CZO}(\ell+\epsilon)$
for all $\ell \in \Z_+$,
where $0<\epsilon<1$
(\cite[p.151, p.154]{Torres}).
\begin{lem}[{\cite[p.152]{Torres}}]\label{4.4}
Let $\sigma \in S_{1,1}^0$.
Then
\[
\sigma(X,D)(x^{\beta})
=(-i)^{|\beta|}\,
\partial_{\xi}^{\beta}(e^{ix\cdot\xi}\sigma(x,\xi)) \bigg|_{\xi=0}
\qquad \text{for all $\beta \in \Z_+^n$}.
\]
\end{lem}
\begin{lem}\label{4.5}
Let $\sigma$ and $\sigma^*$ be defined by {\rm \eqref{(3.19)}}
with $m \le 0$, $0<\delta<1$
and $j_0$ satisfying \eqref{(3.20)} and
$2^{j_0(\delta-1)+2}<2^{-1/2}$,
and {\rm \eqref{(3.1)}}, respectively.
Then the following are true:
\begin{enumerate}
\item
For all $\ell \in \Z_+$,
$\sigma(X,D)$,
$\sigma(X,D)^* \in {\rm CZO}(\ell+\epsilon)$,
where $0<\epsilon<1$,
and for all $\beta \in \Z_+^n$,
$\sigma(X,D)(x^{\beta})=\sigma(X,D)^*(x^{\beta})=0$.
\item
For all $\ell \in \Z_+$,
$\sigma^*(X,D)$,
$\sigma^*(X,D)^* \in {\rm CZO}(\ell+\epsilon)$,
where $0<\epsilon<1$,
and for all $\beta \in \Z_+^n$,
$\sigma^*(X,D)(x^{\beta})=\sigma^*(X,D)^*(x^{\beta})=0$.
\end{enumerate}
\end{lem}
\begin{proof}
Let $m \le 0$.
By Lemma \ref{4.2},
we have
$\sigma, \sigma^* \in S_{1,\delta}^m$
and
$(\partial_{\xi}^{\alpha}\sigma)(x,0)=
(\partial_{\xi}^{\alpha}\sigma^*)(x,0)=0$
for all $\alpha \in \Z_+^n$.
Then,
by $S_{1,\delta}^m \subset S_{1,1}^0$
and the remark below Lemma \ref{4.3},
we get
$\sigma(X,D), \sigma(X,D)^* \in {\rm CZO}(\ell+\epsilon)$
for all $\ell \in \Z_+$,
where $0<\epsilon<1$.
On the other hand,
by Lemma \ref{4.4},
$(\partial_{\xi}^{\alpha}\sigma)(x,0)=0$
for all $\alpha \in \Z_+^n$
gives $\sigma(X,D)(x^{\beta})=0$
for all $\beta \in \Z_+^n$.
By \eqref{(3.2)},
we have
\begin{align*}
\langle \sigma(X,D)^*f,g \rangle
&=\langle \sigma(X,D)g,f \rangle
=(\sigma(X,D)g,\overline{f}) \\
&=(g,\sigma^*(X,D)\overline{f})
=\langle g,\overline{\sigma^*(X,D)\overline{f}} \rangle
=\langle \overline{\sigma^*}(X,-D)f,g \rangle
\end{align*}
for all $f,g \in \calS(\R^n)$.
This implies $\sigma(X,D)^*=\overline{\sigma^*}(X,-D)$.
Hence, by Lemma \ref{4.4}
and $(\partial^{\alpha}\sigma^*)(x,0)=0$ for all $\alpha \in \Z_+^n$,
we see that
\[
\sigma(X,D)^*(x^{\beta})
=(-i)^{|\beta|}\,
\partial_{\xi}^{\beta}(e^{ix\cdot\xi}
\overline{\sigma^*(x,-\xi)}) \bigg|_{\xi=0}=0
\]
for all $\beta \in \Z_+^n$.
Therefore, we obtain Lemma \ref{4.5} (1).

In the same way,
we can prove Lemma \ref{4.5} (2).
Note that $\sigma^*(X,D)^*=\overline{\sigma}(X,-D)$.
\end{proof}
We are now ready to prove Theorem \ref{4.1}.

\medskip
\noindent
{\it Proof of Theorem \ref{4.1}}.
Let $1<p,q <\infty$ and $q\neq 2$.
We only consider the case $q>2$ since
the case $q<2$ can be carried out in a similar way.
By Theorem \ref{3.6},
we see that 
$\sigma$ belongs to $S_{1,\delta}^m$,
but $\sigma(X,D)$ is not bounded on $M^{p,q}(\R^n)$,
where $\sigma$ is defined by \eqref{(3.19)}
with $m>(1/q-1/2)\delta n$,
$0<\delta<1$ and a sufficiently large integer $j_0$.
Since $1/q-1/2<0$,
we may assume $m<0$.
Then, by Lemma \ref{4.5} (1),
we see that
$\sigma(X,D),\sigma(X,D)^* \in {\rm CZO}(\ell+\epsilon)$
for all $\ell \in \Z_+$,
where $0<\epsilon<1$,
and $\sigma(X,D)(x^{\beta})=\sigma(X,D)^*(x^{\beta})=0$
for all $\beta \in \Z_+^n$.
On the other hand,
it is well known that,
if $\delta<1$, then pseudo-differential operators
with symbols in $S_{1,\delta}^0$
are bounded on $L^2(\R^n)$
(\cite[Chapter 2, Theorem 4.1]{Kumano-go},
\cite[Chapter 7, Theorem 2]{Stein}).
It is easy to prove that
the $L^2$-boundedness of $\sigma(X,D)$
gives the weak boundedness property of $\sigma(X,D)$.
Hence,
$\sigma(X,D),\sigma(X,D)^* \in {\rm CZO}(\ell+\epsilon)$
for all $\ell \in \Z_+$,
$\sigma(X,D)(P)=\sigma(X,D)^*(P)=0$
for all polynomials $P$,
and $\sigma(X,D)$ satisfies the weak boundedness property,
but $\sigma(X,D)$ is not bounded on $M^{p,q}(\R^n)$.

\section*{Acknowledgement}
The authors gratefully acknowledge helpful discussions
with Professor Eiichi Nakai and Professor K\^oz\^o Yabuta.


\end{document}